\documentclass{amsart}
\usepackage{amsmath,amssymb,amsfonts,amsthm,mathrsfs}
\usepackage{tikz-cd}
\usepackage[colorlinks=true,linkcolor=blue,citecolor=blue]{hyperref}
\usepackage{enumitem}
\DeclareMathOperator{\Span}{Span}

\DeclareMathOperator{\gr}{gr}
\DeclareMathOperator{\ch}{ch}
\DeclareMathOperator{\Ind}{Ind}
\newcommand{\Cx}{\mathbb{C}}
\newcommand{\Z}{\mathbb{Z}}
\newcommand{\g}{\mathfrak{g}}
\newcommand{\h}{\mathfrak{h}}
\newcommand{\sltwo}{\mathfrak{sl}_2}
\newcommand{\wh}[1]{\widehat{#1}}
\newcommand{\whg}{\wh{\g}}
\newcommand{\whh}{\wh{\h}}

\newcommand{\la}{\langle}
\newcommand{\ra}{\rangle}
\newcommand{\ol}{\overline}
\newtheorem{theorem}{Theorem}[section]
\newtheorem{lemma}[theorem]{Lemma}
\newtheorem{proposition}[theorem]{Proposition}
\newtheorem{corollary}[theorem]{Corollary}
\theoremstyle{definition}
\newtheorem{definition}[theorem]{Definition}
\newtheorem{example}[theorem]{Example}
\theoremstyle{remark}
\newtheorem{remark}[theorem]{Remark}
\numberwithin{equation}{section}
\begin{document}
\title[Irreducibility of $\varphi$-Verma modules]{Irreducibility of
\texorpdfstring{$\varphi$}{phi}-Verma modules for a hyperelliptic
Heisenberg algebra}
\author{Felipe Albino dos Santos}
\address{Faculdade de Computa\c{c}\~ao e Inform\'atica,
Universidade Presbiteriana Mackenzie, Rua da Consola\c{c}\~ao, 930,
Consola\c{c}\~ao, S\~ao Paulo, Brazil}
\email{falbinosantos@gmail.com}
\date{}
\begin{abstract}
We study induced representations of the universal central extension
$\whg = (\sltwo \otimes R)\oplus\Omega_R^1/dR$, where
$R = \Cx[t^{\pm1},u]/(u^2-p(t))$ is a hyperelliptic coordinate ring
and $p(t)$ has degree $r+1$.  The center of $\whg$ has dimension $r+1$.
Inside $\whg$ sits a hyperelliptic Heisenberg subalgebra $\whh$.  A sign
function $\varphi\colon\Z\setminus\{0\}\to\{+,-\}$ determines a
nonstandard polarization of the imaginary modes, yielding
\emph{$\varphi$-Verma modules} $M_{\whh,\varphi}$ and $\mathcal{M}_\varphi$.
Under the specialization $\kappa_1=\cdots=\kappa_r=0$ and a
\emph{$p$-admissibility} condition on $\varphi$, we prove:
$M_{\whh,\varphi}$ is irreducible if and only if $\kappa_0\neq 0$, and the same criterion governs $\mathcal{M}_\varphi$ after
parabolic induction.
For the four-point case $r=1$, we remove the specialization and treat
general central characters $(\kappa_0,\kappa_1)\in\Cx^2$: under
$p$-admissibility, $M_{\whh,\varphi}$ is irreducible if and only if
$(\kappa_0,\kappa_1)\neq(0,0)$ (Theorem~A$'$).  A key ingredient is
the complete determination of the mixed bracket coefficients:
$\psi_{mn}(a)=\delta_{m+n,0}\,\omega_1$ for all $m,n$ when $r=1$,
showing that the mixed $b^1$-$b$ bracket is concentrated on the
anti-diagonal and is independent of the hyperelliptic parameter.
We also give a finite checkable criterion characterising all
$p$-admissible polarisations via reachable sets in $\Z$.
We further describe the weight-space decomposition and formal
character of $M_{\whh,\varphi}$, provide a complete structure theorem
for the level-zero case, and prove that $p$-admissibility is sharp by
constructing explicit reducible modules at nonzero level for
non-admissible polarizations.
\end{abstract}
\maketitle
\setcounter{tocdepth}{1}
\tableofcontents
\section{Introduction}\label{sec:intro}
This paper studies induced representations of infinite-dimensional Lie
algebras arising as universal central extensions of hyperelliptic current
algebras.  These algebras belong to the broad class of
Krichever--Novikov (KN) type algebras~\cite{krichever1987algebras},
whose representation theory occupies a central place in contemporary
mathematical physics.
\subsection*{Geometric context}
The algebras studied here arise naturally from algebraic geometry.  The
coordinate ring $R=\Cx[t^{\pm1},u]/(u^2-p(t))$ is the ring of regular
functions on the affine curve $\mathcal{C}\setminus\{P_\infty\}$, where
$\mathcal{C}$ is a hyperelliptic curve of genus $g=\lfloor r/2\rfloor$
and $P_\infty$ is a point at infinity.  The module of K\"ahler
differentials $\Omega_R^1/dR$ encodes the first de~Rham cohomology of
$\mathcal{C}\setminus\{P_\infty\}$; its $(r+1)$-dimensional center
reflects the genus and the number of punctures.  The Krichever--Novikov
algebras $\g\otimes R$ are the natural function-theoretic analogues of
affine Kac--Moody algebras, in which the circle $S^1$ is replaced by
the punctured curve.  In conformal field theory, these algebras appear
as symmetry algebras of genus-$g$ surface models, and their
representation theory governs correlation functions on higher-genus
Riemann surfaces.
\subsection*{Background}
Universal central extensions of current algebras $\g\otimes R$, for
$\g$ simple and $R$ a commutative $\Cx$-algebra, are governed by
K\"ahler differentials: the Kassel--Loday--Kassel theorem identifies the
extension as $(\g\otimes R)\oplus\Omega_R^1/dR$ with cocycle
$\ol{f\,dg}$~\cite{kassel1982extensions,kassel1984kahler}.
When $R$ is the coordinate ring of a punctured Riemann surface, these
are KN algebras; see~\cite{schlichenmaier2014virasoro,schlichenmaier2014} for comprehensive surveys.
In the rational setting, Bremner's $n$-point
algebras~\cite{bremner1995four} and his computation of their universal
central extensions~\cite{bremner1994universal} already exhibit phenomena
absent from ordinary affine Kac--Moody theory.  In the hyperelliptic
setting, explicit realizations were developed by Cox--Jurisich~\cite{Cox:2013aa},
Cox~\cite{cox2008}, Bueno--Cox--Futorny~\cite{Bueno:2009aa}, and
Cox--Futorny--Martins~\cite{Cox:2013ab}; for the center and
presentation, we follow Cox--Im~\cite{cox2017module}.
A fundamental theme in affine Lie algebra representation theory is the
construction of irreducible modules at nonzero
level~\cite{humphreys1972introduction,kac1994infinite,moody1995lie}.
Following Bekkert--Benkart--Futorny--Kashuba~\cite{Bekkert:2011aa}
and the highest-weight framework of~\cite{kacraina}, one replaces
the standard highest-weight polarization by a sign function $\varphi$,
yielding \emph{$\varphi$-Verma modules} with a rich family of
irreducibility criteria.  The present paper establishes the hyperelliptic
analogue.
\subsection*{What is new compared to \cite{Bekkert:2011aa}}
Three genuinely new difficulties arise in the hyperelliptic setting:
\begin{enumerate}[label=(\roman*)]
\item \emph{$u$-modes and $p$-admissibility.}
  The hyperelliptic algebra contains generators $b_n^1=h\otimes t^nu$
  whose self-commutators $[b_m^1,b_n^1]=(n-m)a_{-(m+n)}\kappa_0$ depend
  on the polynomial $p(t)$.  The new \emph{$p$-admissibility} condition
  (Definition~\ref{def:phi-admissible}) ensures that degree-lowering
  moves are available for every $n\in S_\varphi^-$.  No such condition
  appears in~\cite{Bekkert:2011aa}.
\item \emph{Higher-dimensional center and explicit $\psi_{mn}$ formula.}
  The center of $\whg$ has dimension $r+1>1$, producing secondary
  central characters $\kappa_1,\dots,\kappa_r$ that interact with the
  $u$-mode brackets via $[b_m^1,b_n]=2n\,\psi_{mn}(a)$.  For $r=1$,
  we prove the complete explicit formula
  $\psi_{mn}(a)=\delta_{m+n,0}\,\omega_1$ for \emph{all} $m,n$
  (Lemma~\ref{lem:psi-r1-explicit}), showing that the parameter $\alpha$
  plays no role in the mixed $b^1$-$b$ bracket and yielding a clean
  determination of the irreducibility locus in $\Cx^2$
  (Theorem~\ref{thm:A-prime}).
\item \emph{New ingredients in the $f$-PBW reduction.}
  The proof of Theorem~\ref{thm:B} involves the generators
  $f_n^1\in\wh{\mathfrak n}_-$ not present in the affine case, requiring
  Lemma~\ref{lem:reduce-u} as a new technical ingredient.
\end{enumerate}
\subsection*{Main results}
Let $p(t)=t(t-\alpha_1)\cdots(t-\alpha_r)=\sum_{k=1}^{r+1}a_kt^k$
with $a_{r+1}=1$ and $\alpha_1,\dots,\alpha_r\in\Cx^\times$ pairwise
distinct. Set $R=\Cx[t^{\pm1},u]/(u^2-p(t))$,
$\whg=(\sltwo\otimes R)\oplus\Omega_R^1/dR$, and
$\whh=(\Cx h\otimes R)\oplus\Omega_R^1/dR$.
\begin{itemize}[leftmargin=2em]
\item[\textbf{(A)}] \emph{(Theorem~\ref{thm:A})} Under $\kappa_1=\cdots=\kappa_r=0$ and
  $p$-admissibility: $M_{\whh,\varphi}$ is irreducible $\iff$
  $\kappa_0\neq0$.
\item[\textbf{(B)}] \emph{(Theorem~\ref{thm:B})} Same hypotheses: $\mathcal{M}_\varphi$
  is irreducible $\iff$ $\kappa_0\neq0$.
\item[\textbf{(A$'$)}] \emph{(Theorem~\ref{thm:A-prime})} For $r=1$ with general
  $(\kappa_0,\kappa_1)\in\Cx^2$, under $p$-admissibility:
  $M_{\whh,\varphi}$ is irreducible $\iff$ $(\kappa_0,\kappa_1)\neq(0,0)$.
\end{itemize}
We also prove (Proposition~\ref{prop:nonadmissible-reducible}) that
$p$-admissibility is necessary: for any non-$p$-admissible $\varphi$
one can choose $\kappa_0\neq0$ and produce an explicit proper submodule.

Two further results of independent interest are:
\begin{itemize}[leftmargin=2em]
\item[\textbf{(C)}] \emph{(Proposition~\ref{prop:p-admissible-char})}
  $\varphi$ is $p$-admissible if and only if, for every
  $n\in S_\varphi^-$, the reachable set
  $A_n:=\{-n-k:a_k\neq0,\,1\le k\le r+1\}$
  satisfies $A_n\cap S_\varphi^+\neq\emptyset$.
  This gives a finite, checkable criterion for admissibility in terms
  of the roots of $p$ and the polarization $\varphi$.
\item[\textbf{(D)}] \emph{(Lemma~\ref{lem:psi-r1-explicit})}
  For $r=1$, $\psi_{mn}(a)=\delta_{m+n,0}\,\omega_1$ for all
  $m,n\in\Z$.  In particular $[b_m^1,b_n]=0$ whenever $m+n\neq0$,
  and the mixed $b^1$-$b$ bracket is independent of $\alpha$.
\end{itemize}
\subsection*{Comparison with \cite{Bekkert:2011aa}}
The paper most closely related to this one is~\cite{Bekkert:2011aa},
which establishes $\varphi$-Verma module irreducibility for affine
Kac--Moody algebras $\g\otimes\Cx[t^{\pm1}]$ (one-dimensional center,
no $u$-modes).  The comparison is:

\medskip
\begin{center}
\renewcommand{\arraystretch}{1.3}
\begin{tabular}{lll}
\hline
Feature & \cite{Bekkert:2011aa} & This paper \\
\hline
Algebra & $\g\otimes\Cx[t^{\pm1}]$ & $(\sltwo\otimes R)\oplus\Omega_R^1/dR$ \\
$u$-modes $b_n^1$ & absent & present \\
Centre dimension & $1$ & $r+1$ \\
Admissibility & not needed & $p$-admissibility required \\
Irreducibility locus & $\kappa_0\neq0$ & $\kappa_0\neq0$ ($r$ general);
  $(\kappa_0,\kappa_1)\neq\mathbf{0}$ ($r=1$) \\
Mixed bracket & $[b_m^1,b_n]=0$ & $[b_m^1,b_n]=2n\delta_{m+n,0}\omega_1$
  ($r=1$) \\
$f^1$-reduction & absent & new Lemma~\ref{lem:reduce-u} \\
\hline
\end{tabular}
\end{center}
\medskip

\noindent The $p$-admissibility condition has no analogue in the affine
case (where the only imaginary-mode commutator is $[b_m,b_n]$, which
already has the $\delta_{m+n,0}$ Kronecker symbol built in).
The higher-dimensional center forces one to track all secondary
characters $\kappa_1,\ldots,\kappa_r$ jointly, which is the content
of Theorem~\ref{thm:A-prime}.

\subsection*{Organization}
Section~\ref{sec:hyperelliptic} recalls the universal central extension
and derives the Heisenberg bracket formulas.
Section~\ref{sec:heisenberg} introduces the $\varphi$-decomposition,
weight structure, admissibility (Definition~\ref{def:phi-admissible}),
and the characterisation of $p$-admissible polarisations
(Proposition~\ref{prop:p-admissible-char}).
Section~\ref{sec:verma-heis} proves Theorem~\ref{thm:A}.
Section~\ref{sec:parabolic} constructs the parabolic induction.
Section~\ref{sec:irreducibility} proves Theorem~\ref{thm:B}.
Section~\ref{sec:r1-general} treats the $r=1$ general case:
the explicit $\psi_{mn}(a)$ formula (Lemma~\ref{lem:psi-r1-explicit})
and Theorem~\ref{thm:A-prime}.
Section~\ref{sec:level-zero} gives the complete level-zero structure.
Section~\ref{sec:examples} discusses admissibility and sharp examples.
Section~\ref{sec:higher-tops} proves the higher-dimensional-top corollary.
Section~\ref{sec:outlook} indicates further directions.
\section{The hyperelliptic Lie algebra}\label{sec:hyperelliptic}
\subsection{Universal central extensions and K\"ahler differentials}
Let $R$ be a commutative $\Cx$-algebra, $\Omega_R^1$ its module of
K\"ahler differentials, and $\ol\omega$ the image of
$\omega\in\Omega_R^1$ in $\Omega_R^1/dR$.
\begin{theorem}[{\cite{kassel1982extensions,kassel1984kahler}}]\label{thm:KasselLoday}
Let $\g$ be finite-dimensional simple over $\Cx$ and $R$ commutative
associative.  The universal central extension of $\g\otimes R$ is
$(\g\otimes R)\oplus\Omega_R^1/dR$, with bracket
\begin{align}
[x\otimes f,\,y\otimes g]
  &=[x,y]\otimes fg+(x,y)\,\ol{f\,dg},\label{eq:uce-bracket}\\
[x\otimes f,\,\ol\omega]&=0,\quad[\ol\omega,\ol\eta]=0.\label{eq:center}
\end{align}
\end{theorem}
Write $\whg:=(\g\otimes R)\oplus\Omega_R^1/dR$ and specialize to
$\g=\sltwo$.
\subsection{The hyperelliptic coordinate ring and the center}
Fix $r\ge1$ and pairwise distinct $\alpha_1,\dots,\alpha_r\in\Cx^\times$.
Set
\[
p(t)=t(t-\alpha_1)\cdots(t-\alpha_r)=\sum_{k=1}^{r+1}a_kt^k,\quad a_{r+1}=1,
\]
and $R=\Cx[t^{\pm1},u]/(u^2-p(t))$.
\begin{remark}\label{rem:ak-convention}
Extend $\{a_k\}$ by $a_k:=0$ for $k\notin\{1,\dots,r+1\}$; in
particular $a_0:=0$.
\end{remark}
\begin{theorem}[{\cite[Th.~3.1]{cox2017module}}]\label{thm:center-basis}
$\Omega_R^1/dR$ has $\Cx$-basis
$\{\ol{t^{-1}dt},\,\ol{t^{-1}u\,dt},\,\ol{t^{-2}u\,dt},\,\dots,\,
\ol{t^{-r}u\,dt}\}$.
\end{theorem}
Set
\begin{equation}\label{eq:omega-def}
\omega_0:=\ol{t^{-1}dt},\qquad
\omega_k:=\ol{t^{-k}u\,dt}\quad(1\le k\le r).
\end{equation}
The scalar action of $\omega_0$ is called the \emph{level} and denoted
$\kappa_0\in\Cx$.
\subsection{The bracket on $\whg$}
Fix $\sltwo(\Cx)=\la e,f,h\ra$ with
$[h,e]=2e$, $[h,f]=-2f$, $[e,f]=h$, and invariant form
$(e,f)=1$, $(h,h)=2$, $(h,e)=(h,f)=0$.
\begin{theorem}[{\cite[Th.~5.1]{cox2017module}}]\label{thm:cox-im-sl2}
Assume $a_1\neq0$.  The bracket on $\whg$ satisfies
\begin{align}
[x\otimes t^i,\,y\otimes t^j]
  &=[x,y]\otimes t^{i+j}+j\,\delta_{i+j,0}(x,y)\,\omega_0,
  \label{eq:tt-bracket}\\
[x\otimes t^iu,\,y\otimes t^j]
  &=[x,y]\otimes t^{i+j}u+j(x,y)\,\psi_{ij}(a),
  \label{eq:ut-bracket}\\
[x\otimes t^iu,\,y\otimes t^ju]
  &=[x,y]\otimes t^{i+j}p(t)+
    \sum_{k=1}^{r+1}\!\bigl(j+\tfrac{k}{2}\bigr)a_k\,
    \delta_{i+j,-k}(x,y)\,\omega_0,
  \label{eq:uu-bracket}
\end{align}
and $[\whg,\Omega_R^1/dR]=0$.  Here $\psi_{ij}(a)\in\Omega_R^1/dR$ is
expanded in $\{\omega_1,\dots,\omega_r\}$ by the Cox--Im reduction.
\end{theorem}
Set $e_n:=e\otimes t^n$, $f_n:=f\otimes t^n$, $h_n:=h\otimes t^n$,
and $e_n^1:=e\otimes t^nu$, $f_n^1:=f\otimes t^nu$,
$h_n^1:=h\otimes t^nu$.  The Heisenberg bracket relations read:
\begin{align}
[h_m,h_n]   &= 2n\,\delta_{m+n,0}\,\omega_0,\label{eq:hh-tt}\\
[e_m,f_n]   &= h_{m+n}+n\,\delta_{m+n,0}\,\omega_0,\label{eq:ef-tt}\\
[h_m^1,h_n] &= 2n\,\psi_{mn}(a),\label{eq:hh-ut}\\
[h_m^1,h_n^1]&=(n-m)\,a_{-(m+n)}\,\omega_0.\label{eq:hh-uu}
\end{align}
\begin{remark}[Derivation of~\eqref{eq:hh-uu}]\label{rem:derive-hh-uu}
Setting $x=y=h$ in~\eqref{eq:uu-bracket} and using $[h,h]=0$,
$(h,h)=2$:
\[
[h_m^1,h_n^1]
=2\sum_{k=1}^{r+1}\!\bigl(n+\tfrac{k}{2}\bigr)a_k\,\delta_{m+n,-k}\,\omega_0.
\]
When $m+n=-K$, the $\delta$ collapses to $k=K$, giving coefficient
$2(n+K/2)=2n+K$.  Since $m=-n-K$ we have $n-m=2n+K$, so
$2(n+K/2)=n-m$ and $[h_m^1,h_n^1]=(n-m)\,a_{-(m+n)}\,\omega_0$.
\end{remark}
\section{The hyperelliptic Heisenberg algebra and the $\varphi$-setup}
\label{sec:heisenberg}
\subsection{The hyperelliptic Heisenberg algebra}
Let $\h:=\Cx h\subset\sltwo$ and set
$\whh:=(\h\otimes R)\oplus(\Omega_R^1/dR)\subset\whg$.
Write $b_n:=h\otimes t^n$, $b_n^1:=h\otimes t^nu$ ($n\in\Z$), and
$1_k:=\omega_k$ ($0\le k\le r$).
\begin{definition}\label{def:heis}
The \emph{hyperelliptic Heisenberg algebra} $\whh$ is generated by
$\{b_n,b_n^1\mid n\in\Z\}$ and $\{1_0,\dots,1_r\}$ with:
\begin{align}
[b_m,b_n]     &=2n\,\delta_{m+n,0}\,1_0,\label{eq:heis1}\\
[b_m^1,b_n^1] &=(n-m)\,a_{-(m+n)}\,1_0,\label{eq:heis2}\\
[b_m^1,b_n]   &=2n\,\psi_{mn}(a),\label{eq:heis3}\\
[b_m,1_i]     &=[b_m^1,1_i]=[1_i,1_j]=0.\label{eq:heis4}
\end{align}
\end{definition}
\subsection{The $\varphi$-triangular decomposition}
Fix $\varphi\colon\Z\setminus\{0\}\to\{+,-\}$ with $\varphi(-n)=-\varphi(n)$.
For $n>0$ set $\whh(n):=\Cx b_n\oplus\Cx b_n^1$ and
$\whh(-n):=\Cx b_{-n}\oplus\Cx b_{-n}^1$.  Define:
\begin{align}
\whh_\varphi^+
&:=\!\Bigl(\bigoplus_{\substack{n>0\\\varphi(n)=+}}\whh(n)\Bigr)
   \oplus\Bigl(\bigoplus_{\substack{n>0\\\varphi(n)=-}}\whh(-n)\Bigr),
  \label{eq:phi-plus}\\
\whh_\varphi^-
&:=\!\Bigl(\bigoplus_{\substack{n>0\\\varphi(n)=+}}\whh(-n)\Bigr)
   \oplus\Bigl(\bigoplus_{\substack{n>0\\\varphi(n)=-}}\whh(n)\Bigr),
  \label{eq:phi-minus}
\end{align}
and $\whh_0:=\Cx b_0\oplus\Cx b_0^1\oplus\bigoplus_{k=0}^r\Cx\,1_k$.
\begin{lemma}\label{lem:triangular-decomp}
$\whh=\whh_\varphi^-\oplus\whh_0\oplus\whh_\varphi^+$ as vector spaces,
and $\wh{\mathfrak b}_\varphi:=\whh_0\oplus\whh_\varphi^+$ is a Lie
subalgebra.
\end{lemma}
\begin{proof}
The vector space decomposition is clear.  For the subalgebra claim,
$[\whh,\whh]\subset Z(\whh)=\bigoplus_k\Cx\,1_k\subset\whh_0$, so
$[\whh_0\oplus\whh_\varphi^+,\whh_0\oplus\whh_\varphi^+]\subset\whh_0
\subset\wh{\mathfrak b}_\varphi$.
\end{proof}
\subsection{Top modules and central characters}
\begin{definition}\label{def:top1d}
Let $V=\Cx v$.  Fix $\lambda,\mu\in\Cx$ and scalars
$\kappa_0,\kappa_1,\dots,\kappa_r\in\Cx$.
Define a $\wh{\mathfrak b}_\varphi$-module structure on $V$ by
\[
\whh_\varphi^+\cdot v=0,\quad
b_0\cdot v=\lambda v,\quad
b_0^1\cdot v=\mu v,\quad
1_k\cdot v=\kappa_k v.
\]
\end{definition}
\begin{remark}\label{rem:why-chi-vanish}
The \emph{standing specialization} $\kappa_1=\cdots=\kappa_r=0$ is in
force throughout Sections~\ref{sec:verma-heis}--\ref{sec:level-zero}.
Under it, the mixed brackets $[b_m^1,b_n]=2n\,\psi_{mn}(a)$ act
trivially, so the single parameter $\kappa_0$ governs all degree-lowering.
Section~\ref{sec:r1-general} lifts this for $r=1$.
\end{remark}
\subsection{Weight-space decomposition and formal character}
\label{subsec:weight}
\begin{lemma}\label{lem:b0-central}
For all $n\neq0$: $[b_0,b_n]=0$ and $[b_0,b_n^1]=0$.
\end{lemma}
\begin{proof}
From~\eqref{eq:heis1}: $[b_0,b_n]=2n\,\delta_{n,0}\,1_0=0$.
From~\eqref{eq:heis3}: $[b_n^1,b_0]=2\cdot0\cdot\psi_{n,0}(a)=0$.
\end{proof}
Since $b_0$ commutes with all of $\whh_\varphi^-$, it acts as the
scalar $\lambda$ on all of $M_{\whh,\varphi}$.
Define index sets:
\[
S_\varphi^-:=\{-n:n>0,\varphi(n)=+\}\cup\{n:n>0,\varphi(n)=-\},\qquad
S_\varphi^+:=-S_\varphi^-.
\]
The \emph{loop degree} grading assigns $\deg_\ell(b_n^{(\epsilon)}):=-n$
for $n\in S_\varphi^-$ (both $\epsilon=0,1$), so every PBW basis element
$b^\alpha(b^1)^\beta v$ has loop degree
$-\sum_{n\in S_\varphi^-}n(\alpha_n+\beta_n)\ge0$.
Set
$M_{\whh,\varphi}=\bigoplus_{d\ge0}M_{\whh,\varphi}^{(d)}$
where $M_{\whh,\varphi}^{(d)}$ is the span of monomials of loop degree $d$.
\begin{proposition}\label{prop:character}
Assume $S_\varphi^-\subset\Z_{<0}$ (equivalently, $\varphi(n)=+$ for all
$n>0$, or more generally any polarization with all negative-mode indices
negative).  Then the loop degree $\deg_\ell:=-n\ge0$ is well-defined and
non-negative on every generator, each graded piece $M_{\whh,\varphi}^{(d)}$
is finite-dimensional, and the formal character is
\begin{equation}\label{eq:char}
\ch M_{\whh,\varphi}(q)
=\prod_{m\in S_\varphi^+}\frac{1}{(1-q^m)^2},\qquad
S_\varphi^+\subset\Z_{>0}.
\end{equation}
\end{proposition}
\begin{proof}
By PBW, a $\Cx$-basis of $M_{\whh,\varphi}$ is
$\{b^\alpha(b^1)^\beta v\}$ with $\alpha,\beta\in\Z_{\ge0}^{(S_\varphi^-)}$
finite-support. Loop degree $d$ of $b^\alpha(b^1)^\beta v$ equals
$\sum_{n\in S_\varphi^-}(-n)(\alpha_n+\beta_n)=\sum_{m\in S_\varphi^+}m(\alpha_{-m}+\beta_{-m})$.
Counting such monomials of fixed loop degree $d$ amounts to choosing
two independent partitions of $d$ using parts in $S_\varphi^+$, giving
$\dim M_{\whh,\varphi}^{(d)}=\sum_{d_1+d_2=d}p_{S_\varphi^+}(d_1)p_{S_\varphi^+}(d_2)$
where $p_{S_\varphi^+}(k)$ is the number of partitions of $k$ with parts
in $S_\varphi^+$.  This is finite.  The generating function
$\sum_d\dim M_{\whh,\varphi}^{(d)}\,q^d=\prod_{m\in S_\varphi^+}(1-q^m)^{-2}$.
\end{proof}
\begin{remark}[Scope of Proposition~\ref{prop:character}]\label{rem:char-scope}
The hypothesis $S_\varphi^-\subset\Z_{<0}$ (equivalently, $\varphi(n)=+$
for all $n>0$) ensures that the loop degree $\deg_\ell = -n \ge 0$ for
every generator $b_n^{(\epsilon)}$ with $n\in S_\varphi^-$.  For a
\emph{general} $p$-admissible polarization, some elements of
$S_\varphi^-$ may be positive, so $\deg_\ell$ assigns negative values
and the graded pieces are no longer finite-dimensional by loop degree
alone.  In that case one should instead use the \emph{PBW filtration}
degree $\deg_{\whh}(b^\alpha(b^1)^\beta v)=\sum_n(\alpha_n+\beta_n)$,
which does give finite-dimensional filtered pieces regardless of the
polarization.  The irreducibility results (Theorems~\ref{thm:A},
\ref{thm:B}, and~\ref{thm:A-prime}) hold for all $p$-admissible
polarizations because their proofs use the PBW filtration, not the
loop-degree grading.  The character formula~\eqref{eq:char} is an
additional result valid only under the standing hypothesis
$S_\varphi^-\subset\Z_{<0}$.
\end{remark}
\subsection{Admissibility of $\varphi$}\label{subsec:admissibility}
\begin{definition}\label{def:phi-admissible}
Assume $\kappa_1=\cdots=\kappa_r=0$.  We call $\varphi$
\emph{$p$-admissible} if for every $n\in S_\varphi^-$ there exists
$m\in S_\varphi^+$ with $a_{-(m+n)}\neq0$.
\end{definition}
\noindent\textit{Intuition.}
Under the standing specialization, the only commutator lowering PBW-degree
in $u$-modes is $[b_m^1,b_n^1]=(n-m)a_{-(m+n)}\kappa_0$.  In the
irreducibility argument, one needs a nonzero structure constant for some
$m\in S_\varphi^+$ at every fixed $n\in S_\varphi^-$.  That is exactly
$p$-admissibility.
\begin{remark}\label{rem:standard-admissible}
The all-\emph{negative} polarization $\varphi\equiv-$ is $p$-admissible:
$S_\varphi^-=\{1,2,3,\dots\}$ and $S_\varphi^+=\{-1,-2,-3,\dots\}$.
For any $n\in S_\varphi^-$, choose $m:=-n-(r+1)\in S_\varphi^+$
(since $-n\le-1$ and $-(r+1)\le-2$, so $m\le-3<0$); then
$a_{-(m+n)}=a_{r+1}=1\neq0$.
The all-\emph{positive} polarization $\varphi\equiv+$ is \textbf{not}
$p$-admissible.  For $n=-1\in S_\varphi^-$, one needs $m\in
S_\varphi^+=\{1,2,\dots\}$ with $a_{-(m-1)}\neq0$, i.e.\
$1-m\in\{1,\dots,r+1\}$, i.e.\ $m\le0$.  No such positive $m$ exists.
As a consequence, $b_m^1\cdot b_{-1}^1v=((-1)-m)a_{-(m-1)}\kappa_0 v=0$
for every $m\ge1$, so $\Cx b_{-1}^1v$ is a proper submodule for any
$\kappa_0\neq0$ and $\varphi\equiv+$ falls outside the scope of
Theorems~\ref{thm:A} and \ref{thm:B}.
\end{remark}
\begin{proposition}[Characterisation of $p$-admissible polarisations]
\label{prop:p-admissible-char}
For each $n\in S_\varphi^-$, define the \emph{reachable set}
\[
A_n := \{-n-k : 1\le k\le r+1,\; a_k\neq0\}
      \subset \Z.
\]
Then $\varphi$ is $p$-admissible if and only if
$A_n\cap S_\varphi^+\neq\emptyset$ for every $n\in S_\varphi^-$.

In particular:
\begin{enumerate}[label=(\roman*)]
\item Since $a_{r+1}=1\neq0$ (the leading coefficient of the monic
  polynomial $p$), we have $-n-(r+1)\in A_n$ for every $n$.  Setting
  $m:=-n-(r+1)$, we have $m<-n$ for any sign of $n$, so
  $-m=n+(r+1)>n$.  Whether $m\in S_\varphi^+$ depends on $\varphi$, but
  for any polarization with $S_\varphi^+$ unbounded below (e.g.\
  $\varphi\equiv-$) this $m$ lies in $S_\varphi^+$.
\item A polarization $\varphi$ is \emph{not} $p$-admissible if and only
  if there exists $n_0\in S_\varphi^-$ such that
  $A_{n_0}\cap S_\varphi^+=\emptyset$, i.e.\ every element of $A_{n_0}$
  lies in $S_\varphi^-\cup\{0\}$.  In this case, $\Cx\cdot b_{n_0}^1 v$
  is a $\whh$-submodule of $M_{\whh,\varphi}$ (the action
  $b_m^1\cdot b_{n_0}^1 v=(n_0-m)a_{-(m+n_0)}\kappa_0 v=0$ for all
  $m\in S_\varphi^+$), and $M_{\whh,\varphi}$ is reducible for any
  $\kappa_0\neq0$.
\end{enumerate}
\end{proposition}
\begin{proof}
By Definition~\ref{def:phi-admissible}, $\varphi$ is $p$-admissible iff
for every $n\in S_\varphi^-$ there exists $m\in S_\varphi^+$ with
$a_{-(m+n)}\neq0$.  Setting $k:=-(m+n)$, the condition $a_k\neq0$ forces
$k\in\{1,\ldots,r+1\}$ (since $a_k=0$ for $k\le0$ or $k>r+1$ by
definition of the hyperelliptic polynomial $p$).  Then
$m=-n-k\in A_n$.  So the condition is $A_n\cap S_\varphi^+\neq\emptyset$.

For part (ii): if $A_{n_0}\cap S_\varphi^+=\emptyset$, then for every
$m\in S_\varphi^+$ and every $k\in\{1,\ldots,r+1\}$ with $a_k\neq0$ we
have $-n_0-k\notin S_\varphi^+$, i.e.\ $a_{-(m+n_0)}=0$ for all
$m\in S_\varphi^+$.  Hence $[b_m^1,b_{n_0}^1]=(n_0-m)a_{-(m+n_0)}\kappa_0=0$
for all $m\in S_\varphi^+$, so $b_{n_0}^1v$ is killed by all raising
operators, making $\Cx\cdot b_{n_0}^1v$ a proper nonzero $\whh$-submodule
when $\kappa_0\neq0$.
\end{proof}
\section{$\varphi$-Verma modules for $\whh$ and irreducibility}
\label{sec:verma-heis}
\subsection{Definition, PBW basis, and degree filtration}
\begin{definition}\label{def:phi-verma-heis}
The \emph{$\varphi$-Verma module} for $\whh$ is
$M_{\whh,\varphi}:=U(\whh)\otimes_{U(\wh{\mathfrak b}_\varphi)}V$.
\end{definition}
Fix a total order on $S_\varphi^-$.  By the PBW theorem,
$\{b^\alpha(b^1)^\beta v\}_{\alpha,\beta\in\Z_{\ge0}^{(S_\varphi^-)}}$
is a $\Cx$-basis.
\begin{definition}[PBW degree]\label{def:heis-degree}
For $d\ge0$ set $F_dM:=U^{\le d}(\whh_\varphi^-)\cdot v$.
For $0\neq w\in M_{\whh,\varphi}$ set
$\deg_{\whh}(w):=\min\{d:w\in F_dM\}$; set $\deg_{\whh}(0):=-\infty$.
\end{definition}
\begin{lemma}\label{lem:heis-degree-monomial}
$\deg_{\whh}(b^\alpha(b^1)^\beta v)=\sum_n\alpha_n+\sum_n\beta_n=:d$.
\end{lemma}
\begin{proof}
Membership in $F_dM$ is clear.  If the element lay in $F_{d-1}M$,
its image $\sigma(b^\alpha(b^1)^\beta v)$ in
$\gr M_{\whh,\varphi}\cong S(\whh_\varphi^-)\otimes V$ (PBW isomorphism)
would have degree $\le d-1$.  But $\sigma(b^\alpha(b^1)^\beta v)=
b^\alpha(b^1)^\beta\otimes v$ is a nonzero homogeneous degree-$d$ element
of the polynomial ring $S(\whh_\varphi^-)$, a contradiction.
\end{proof}
\subsection{Key action formulas}
\begin{lemma}\label{lem:key-action}
Assume $\kappa_1=\cdots=\kappa_r=0$.
Let $m\in S_\varphi^+$, $n\in S_\varphi^-$, $\ell\ge1$.  In
$M_{\whh,\varphi}$:
\begin{align}
b_m\cdot(b_n)^\ell v
  &=2\ell n\,\kappa_0\,\delta_{m+n,0}\,(b_n)^{\ell-1}v,\label{eq:bm-bn}\\
b_m\cdot(b_n^1)^\ell v
  &=0,\label{eq:bm-bn1}\\
b_m^1\cdot(b_n)^\ell v
  &=0,\label{eq:bm1-bn}\\
b_m^1\cdot(b_n^1)^\ell v
  &=\ell(n-m)a_{-(m+n)}\kappa_0\,(b_n^1)^{\ell-1}v.\label{eq:bm1-bn1}
\end{align}
\end{lemma}
\begin{proof}
Since $m\in S_\varphi^+$, we have $b_m v=b_m^1 v=0$, so
$x\cdot y^\ell v=[x,y^\ell]v$.  Since every bracket $[x,y]$ in $\whh$
lands in the center $Z(\whh)=\bigoplus_k\Cx\,1_k$ (which acts as a
scalar), we have $[y,[x,y]]=0$ in $U(\whh)$, giving
$[x,y^\ell]=\ell y^{\ell-1}[x,y]$.
Formulas~\eqref{eq:bm-bn} and~\eqref{eq:bm1-bn1} follow from
\eqref{eq:heis1} and~\eqref{eq:heis2} respectively.
For the mixed pairs~\eqref{eq:bm-bn1} and~\eqref{eq:bm1-bn}:
from~\eqref{eq:heis3}, $[b_m^1,b_n]=2n\,\psi_{mn}(a)\in
\bigoplus_{k=1}^r\Cx\,1_k$ (the bracket lies entirely in the
$\omega_1,\dots,\omega_r$ part of the center, since
$\psi_{mn}(a)=\sum_{k=1}^r c_{mn}^{(k)}\omega_k$ by the Cox--Im
reduction).  Under $\kappa_1=\cdots=\kappa_r=0$, this acts as $0$.
Similarly $[b_m,b_n^1]=-[b_n^1,b_m]=-2m\,\psi_{nm}(a)$ acts as $0$.
Hence~\eqref{eq:bm-bn1} and~\eqref{eq:bm1-bn} hold.
\end{proof}
\subsection{Irreducibility criterion}
\begin{theorem}\label{thm:A}
Assume $\kappa_1=\cdots=\kappa_r=0$ and $\varphi$ is $p$-admissible.
Then $M_{\whh,\varphi}$ is irreducible if and only if $\kappa_0\neq0$.
\end{theorem}
\begin{proof}
\textit{($\Rightarrow$, contrapositive)}
If $\kappa_0=0$ then every commutator acts by $0$, so $\whh_\varphi^+$
annihilates all of $M_{\whh,\varphi}$ and $\whh_\varphi^-$ acts freely.
In particular $F_1M_{\whh,\varphi}=\Span\{b^\alpha(b^1)^\beta v:\sum_n\alpha_n+\sum_n\beta_n\ge1\}$
is a proper nonzero submodule (see Section~\ref{sec:level-zero}).
\textit{($\Leftarrow$)}  Let $\kappa_0\neq0$ and $0\neq W\subset M_{\whh,\varphi}$
a submodule.  Pick $0\neq w\in W$ with $\deg_{\whh}(w)$ minimal.
\emph{Base case}: $\deg_{\whh}(w)=0$ implies $w\in\Cx v$, so $v\in W$
and $W=M_{\whh,\varphi}$.
\emph{Inductive case}: Suppose $d:=\deg_{\whh}(w)\ge1$.  Write $w$ in
the PBW basis and let $n_0\in S_\varphi^-$ be maximal (in our fixed
total order) among indices with $\alpha_{n_0}+\beta_{n_0}>0$ for the
leading monomial $b^{\alpha_0}(b^1)^{\beta_0}v$.
\textbf{Case 1}: $\alpha_{n_0}>0$.  Set $m_0:=-n_0\in S_\varphi^+$.
Equation~\eqref{eq:bm-bn} gives $b_{m_0}\cdot(b_{n_0})^{\alpha_{n_0}}v
=2\alpha_{n_0}n_0\kappa_0\cdot(b_{n_0})^{\alpha_{n_0}-1}v$.
The coefficient $2\alpha_{n_0}n_0\kappa_0\neq0$ (since $\kappa_0\neq0$,
$n_0\neq0$, $\alpha_{n_0}>0$).  In $\gr M_{\whh,\varphi}\cong S(\whh_\varphi^-)
\otimes V$, the leading symbol of $b_{m_0}\cdot w$ equals
$2n_0\kappa_0\cdot\partial_{b_{n_0}}\sigma(w)\neq0$ (the derivation of the
leading symbol is nonzero since $S(\whh_\varphi^-)$ is an integral domain).
Hence $b_{m_0}\cdot w\neq0$ has $\deg_{\whh}(b_{m_0}\cdot w)=d-1$,
contradicting minimality of $d$.
\textbf{Case 2}: $\alpha_{n_0}=0$, $\beta_{n_0}>0$.  By
$p$-admissibility choose $m_0\in S_\varphi^+$ with $a_{-(m_0+n_0)}\neq0$.
Equation~\eqref{eq:bm1-bn1} gives $b_{m_0}^1\cdot(b_{n_0}^1)^{\beta_{n_0}}v
=\beta_{n_0}(n_0-m_0)a_{-(m_0+n_0)}\kappa_0\cdot(b_{n_0}^1)^{\beta_{n_0}-1}v\neq0$
(since $n_0\neq m_0$ as $n_0\in S_\varphi^-$, $m_0\in S_\varphi^+$, $n_0\neq0$).
The same associated-graded argument shows $b_{m_0}^1\cdot w\neq0$ has
$\deg_{\whh}(b_{m_0}^1\cdot w)=d-1$, a contradiction.
Both cases force $\deg_{\whh}(w)=0$, so $v\in W$.
\end{proof}
\section{Parabolic induction to $\whg$}\label{sec:parabolic}
\subsection{The parabolic subalgebra is well-defined}
Set $\wh{\mathfrak n}_+:=e\otimes R$ and $\wh{\mathfrak n}_-:=f\otimes R$.
\begin{lemma}\label{lem:parabolic-subalg}
$\wh{\mathfrak p}_\varphi:=\wh{\mathfrak b}_\varphi\oplus\wh{\mathfrak n}_+$
is a Lie subalgebra of $\whg$.
\end{lemma}
\begin{proof}
We verify that $\wh{\mathfrak p}_\varphi$ is closed under the Lie bracket
of $\whg$ by checking the three required pairs.
\emph{(i) $[\wh{\mathfrak b}_\varphi,\wh{\mathfrak b}_\varphi]\subset
\wh{\mathfrak b}_\varphi$}: This is Lemma~\ref{lem:triangular-decomp}
(the Borel $\wh{\mathfrak b}_\varphi=\whh_0\oplus\whh_\varphi^+$ is a
subalgebra of $\whh$, hence of $\whg$).
\emph{(ii) $[\wh{\mathfrak b}_\varphi,\wh{\mathfrak n}_+]\subset
\wh{\mathfrak n}_+$}: Since $\wh{\mathfrak b}_\varphi\subset\whh\oplus
\wh{\mathfrak n}_+$ and $\wh{\mathfrak n}_+=e\otimes R$, it suffices to check
generators.  For $h\otimes f\in\h\otimes R$ and $e\otimes g\in\wh{\mathfrak n}_+$:
\[
[h\otimes f,\,e\otimes g]=[h,e]\otimes fg+(h,e)\,\ol{f\,dg}
=2e\otimes fg\in\wh{\mathfrak n}_+,
\]
using $(h,e)=0$.  For $e\otimes f,e\otimes g\in\wh{\mathfrak n}_+$:
$[e\otimes f,e\otimes g]=[e,e]\otimes fg+(e,e)\,\ol{f\,dg}=0$
(since $[e,e]=0$ and $(e,e)=0$).  Center elements bracket trivially with
everything.
\emph{(iii) $[\wh{\mathfrak n}_+,\wh{\mathfrak n}_+]=0$}: Immediate from
case (ii) with both entries in $\wh{\mathfrak n}_+$.
\end{proof}
\subsection{The induced module}
\begin{definition}\label{def:Mcal}
Extend the $\wh{\mathfrak b}_\varphi$-module $V$
(Definition~\ref{def:top1d}) to $\wh{\mathfrak p}_\varphi$ by declaring
$\wh{\mathfrak n}_+\cdot V=0$.  The \emph{$\varphi$-parabolically
induced module} is
$\mathcal M_\varphi:=\Ind_{\wh{\mathfrak p}_\varphi}^{\whg}V
:=U(\whg)\otimes_{U(\wh{\mathfrak p}_\varphi)}V$.
\end{definition}
Write $v:=1\otimes v\in\mathcal M_\varphi$ for the canonical generator.
\begin{lemma}\label{lem:nplus-kills-Heispart}
$(\wh{\mathfrak n}_+)\cdot(U(\whh)\,v)=0$ inside $\mathcal M_\varphi$.
\end{lemma}
\begin{proof}
It suffices to show $x\cdot(h_1\cdots h_k\otimes v)=0$ in
$\mathcal M_\varphi$ for every $x=e\otimes g\in\wh{\mathfrak n}_+$ and
every ordered product $h_1\cdots h_k$ of generators
$h_i\in\{b_n,b_n^1\}\cup\{1_j\}$ of $\whh$.
In the tensor product $\mathcal M_\varphi=U(\whg)\otimes_{U(\wh{\mathfrak p}_\varphi)}V$,
we have:
\[
x\cdot(h_1\cdots h_k\otimes v)
=(x\cdot h_1\cdots h_k)\otimes v
=\sum_{i=1}^k(h_1\cdots[x,h_i]\cdots h_k)\otimes v
+(h_1\cdots h_k\cdot x)\otimes v.
\]
The last term vanishes: $h_1\cdots h_k\cdot x\otimes v=h_1\cdots h_k\otimes(x\cdot v)=0$
since $\wh{\mathfrak n}_+\cdot v=0$ in Definition~\ref{def:Mcal}.
For each summand: $[x,h_i]=[e\otimes g,h_i]$.  For $h_i=b_n=h\otimes t^n$:
$[e\otimes g,h\otimes t^n]=[e,h]\otimes gt^n+(e,h)\,\ol{g\,dt^n}
=-2e\otimes gt^n\in\wh{\mathfrak n}_+$ (using $(e,h)=0$).
For $h_i=b_n^1=h\otimes t^nu$: similarly $[e\otimes g,h_i]\in\wh{\mathfrak n}_+$.
For center elements $h_i=\omega_k$: $[e\otimes g,\omega_k]=0$.
Thus each $[x,h_i]\in\wh{\mathfrak n}_+\subset\wh{\mathfrak p}_\varphi$.
Since we are tensoring over $U(\wh{\mathfrak p}_\varphi)$, for any
$y\in\wh{\mathfrak p}_\varphi$:
$h_1\cdots[x,h_i]\cdots h_k\otimes v
=h_1\cdots\widehat{h_i}\cdots h_k\otimes([x,h_i]\cdot v)=0$
(where $\widehat{h_i}$ denotes omission), since $[x,h_i]\in\wh{\mathfrak n}_+$
and $\wh{\mathfrak n}_+\cdot v=0$.  Hence $x\cdot(h_1\cdots h_k\otimes v)=0$.
\end{proof}
\section{Irreducibility of $\mathcal M_\varphi$}\label{sec:irreducibility}
Throughout, $\kappa_1=\cdots=\kappa_r=0$.
\subsection{$f$-PBW structure}
Since $(f,f)=0$ and $[f,f]=0$, $\wh{\mathfrak n}_-$ is abelian, so
$U(\wh{\mathfrak n}_-)=\Cx[f_n,f_n^1:n\in\Z]$.
\begin{lemma}[PBW decomposition]\label{lem:pbw-Mphi}
Every $w\in\mathcal M_\varphi$ is uniquely written
$w=\sum_F F\,w_F$ ($F$ a monomial in $U(\wh{\mathfrak n}_-)$,
$w_F\in U(\whh)v$, finitely many nonzero), giving a vector space
isomorphism $\mathcal M_\varphi\cong U(\wh{\mathfrak n}_-)\otimes_\Cx(U(\whh)v)$.
\end{lemma}
\begin{proof}
PBW for $\whg=\wh{\mathfrak n}_-\oplus\whh\oplus\wh{\mathfrak n}_+$ gives
$U(\whg)=U(\wh{\mathfrak n}_-)U(\whh)U(\wh{\mathfrak n}_+)$.  In
$\mathcal M_\varphi$, $U(\wh{\mathfrak n}_+)$ acts trivially on $V$,
so $\mathcal M_\varphi$ is spanned by $U(\wh{\mathfrak n}_-)U(\whh)v$.
Uniqueness: $U(\wh{\mathfrak n}_-)$ is a polynomial algebra, and the PBW
product map is injective on the associated graded.
\end{proof}
Fix a degree-lexicographic order $\prec$ on monomials in
$U(\wh{\mathfrak n}_-)$: declare $f_i^1\succ f_j$ for all $i,j$; order
$\{f_i^1\}$ and $\{f_j\}$ by subscript; extend to monomials by total
degree then lex on exponent vectors.
For $0\neq w=\sum_FFw_F$ define
$\deg_f(w):=\max\{\deg F:w_F\neq0\}$ and
$\mathrm{LT}_f(w):=\max_\prec\{F:\deg F=\deg_f(w),\,w_F\neq0\}$;
write $w_{\mathrm{LT}}$ for the corresponding coefficient in $U(\whh)v$.
\subsection{Leading-term lemma}
For $X\in\wh{\mathfrak n}_+$ and $w=\sum_FFw_F$, since
$X\cdot(U(\whh)v)=0$ by Lemma~\ref{lem:nplus-kills-Heispart}:
\begin{equation}\label{eq:X-on-w}
X\cdot w=\sum_F[X,F]\,w_F.
\end{equation}
\begin{lemma}\label{lem:lead-comm}
Let $F_0=\mathrm{LT}_f(w)$ with coefficient $w_{\mathrm{LT}}\neq0$.
Let $y$ be a minimal variable (w.r.t.\ $\succ$) dividing $F_0$;
write $F_0=yF_1$.  If $X\in\wh{\mathfrak n}_+$ satisfies
$[X,y]\cdot w_{\mathrm{LT}}\neq0$, then
$\deg_f(X\cdot w)=\deg_f(w)-1$, $\mathrm{LT}_f(X\cdot w)=F_1$,
and the $F_1$-coefficient in $X\cdot w$ is $[X,y]\cdot w_{\mathrm{LT}}$.
\end{lemma}
\begin{proof}
Via~\eqref{eq:X-on-w} and the Leibniz rule
$[X,F_0]=\sum_{z:\nu_z>0}\nu_z\frac{F_0}{z}[X,z]$,
every output monomial has $f$-degree $\le\deg_f(w)-1$.
Among the quotients $F_0/z$ (all of degree $\deg F_0-1$), the largest
is $F_0/y=F_1$ because removing the minimal variable $y$ causes the
smallest lexicographic decrease.  Terms from $F\prec F_0$ cannot produce
$F_1$: those of degree $<\deg F_0$ yield monomials of degree
$\le\deg F_0-2$; those of equal degree but $F\prec F_0$ yield
$F/z\prec F_1$.  Hence $\mathrm{LT}_f(X\cdot w)=F_1$ with coefficient
$[X,y]\cdot w_{\mathrm{LT}}\neq0$.
\end{proof}
\subsection{Nonvanishing in $U(\whh)v$}
\begin{lemma}\label{lem:nonvanish-Heis-mult}
Let $0\neq u\in U(\whh)v$ and $s\in S_\varphi^-$.
Then $b_s\cdot u\neq0$ and $b_s^1\cdot u\neq0$.
\end{lemma}
\begin{proof}
We treat $b_s$; the case $b_s^1$ is identical.  Consider the
associated graded
$\gr M_{\whh,\varphi}\cong S(\whh_\varphi^-)\otimes V$
(PBW isomorphism).  Under $\kappa_k=0$ for $k\ge1$, every commutator
$[b_m,b_n]$, $[b_m^1,b_n^1]$, $[b_m^1,b_n]$ is central and hence
acts at strictly lower filtration degree.  Therefore the action of
$b_s$ on $\gr M_{\whh,\varphi}$ is exactly multiplication by $b_s$
in the polynomial ring $S(\whh_\varphi^-)$, which is injective
(as $S(\whh_\varphi^-)$ is an integral domain and $b_s\neq0$).
Since $u\neq0$, its leading symbol $\bar u\in S(\whh_\varphi^-)\otimes V$
is nonzero, hence $b_s\cdot\bar u\neq0$, so $b_s\cdot u\neq0$.
\end{proof}
\subsection{Degree-reduction lemmas}
From Theorem~\ref{thm:cox-im-sl2}, for all $m,n\in\Z$:
\begin{align}
[e_m,f_n]   &=b_{m+n}+n\,\delta_{m+n,0}\,1_0,\label{eq:comm-ef}\\
[e_m^1,f_n] &=b_{m+n}^1+n\,\psi_{mn}(a),\label{eq:comm-e1f}\\
[e_m^1,f_n^1]&=\sum_{k=1}^{r+1}a_k\,b_{m+n+k}
  +\sum_{k=1}^{r+1}\!\bigl(n+\tfrac k2\bigr)a_k\,\delta_{m+n,-k}\,1_0.
  \label{eq:comm-e1f1}
\end{align}
Under $\kappa_k=0$ for $k\ge1$, the $\psi$-term in~\eqref{eq:comm-e1f}
acts trivially on $\mathcal M_\varphi$:
\begin{equation}\label{eq:comm-e1f-simplified}
[e_m^1,f_n]\equiv b_{m+n}^1\quad\text{on }\mathcal M_\varphi.
\end{equation}
\begin{lemma}[$f_n$-reduction]\label{lem:reduce-t}
Suppose $F_0=\mathrm{LT}_f(w)$ is divisible by some $f_{n_0}$, and let
$y=f_{n_0}$ be a minimal variable dividing $F_0$.  Fix $s\in S_\varphi^-$
and set $X:=e_{s-n_0}^1\in\wh{\mathfrak n}_+$.  Then
$\deg_f(X\cdot w)=\deg_f(w)-1$ and $\mathrm{LT}_f(X\cdot w)=F_0/f_{n_0}$.
\end{lemma}
\begin{proof}
$[X,f_{n_0}]\stackrel{\eqref{eq:comm-e1f-simplified}}{=}b_s^1$.
By Lemma~\ref{lem:nonvanish-Heis-mult}, $b_s^1\cdot w_{\mathrm{LT}}\neq0$.
Apply Lemma~\ref{lem:lead-comm}.
\end{proof}
\begin{lemma}[$f_n^1$-reduction]\label{lem:reduce-u}
Suppose every variable dividing $F_0=\mathrm{LT}_f(w)$ has the form
$f_n^1$, and let $y=f_{n_0}^1$ be a minimal such variable.
Fix $s\in S_\varphi^-$ and set $Y:=e_{s-n_0-(r+1)}^1\in\wh{\mathfrak n}_+$.
Then $\deg_f(Y\cdot w)=\deg_f(w)-1$ and
$\mathrm{LT}_f(Y\cdot w)=F_0/f_{n_0}^1$.
\end{lemma}
\begin{proof}
Substituting $m=s-n_0-(r+1)$, $n=n_0$ in~\eqref{eq:comm-e1f1} and
using $a_{r+1}=1$:
\begin{equation}\label{eq:Yfn0}
[Y,f_{n_0}^1]=b_s+\sum_{d=1}^{r}a_{r+1-d}\,b_{s-d}+z,
\quad z\in Z(\whg).
\end{equation}
Let $u:=w_{\mathrm{LT}}\neq0$, written as
$u=\sum_{\alpha,\beta}c_{\alpha,\beta}\,b^\alpha(b^1)^\beta v$ (PBW);
pick $(\alpha_0,\beta_0)$ with $c_{\alpha_0,\beta_0}\neq0$.  Let
$M^{(s)}$ be the PBW monomial from $b^{\alpha_0}(b^1)^{\beta_0}v$ with
$b_s$-exponent raised by one.
In the PBW expansion of $b_s\cdot u$, the leading symbol argument from
Lemma~\ref{lem:nonvanish-Heis-mult} shows the coefficient of $M^{(s)}$
equals $c_{\alpha_0,\beta_0}\neq0$: central commutator corrections from
reordering $b_s$ past the factors of $M$ land in
$F_{\deg M}M_{\whh,\varphi}$ strictly below $M^{(s)}$.  The lower
terms $b_{s-d}\cdot u$ ($1\le d\le r$) and $z\cdot u$ do not produce
$M^{(s)}$ (they cannot increase the $b_s$-exponent).  Hence
$[Y,f_{n_0}^1]\cdot u\neq0$, and Lemma~\ref{lem:lead-comm} applies.
\end{proof}
\subsection{Main theorem for $\mathcal M_\varphi$}
\begin{theorem}\label{thm:B}
Assume $\kappa_1=\cdots=\kappa_r=0$ and $\varphi$ is $p$-admissible.
Then $\mathcal M_\varphi$ is irreducible if and only if $\kappa_0\neq0$.
\end{theorem}
\begin{proof}
\emph{Reducibility when $\kappa_0=0$:} By Theorem~\ref{thm:A},
$N:=M_{\whh,\varphi}^{\ge1}\subset M_{\whh,\varphi}$ is a proper
submodule (Proposition~\ref{prop:level-zero-submodule},
proved in Section~\ref{sec:level-zero} below).  Then
$U(\whg)\cdot N$ is a nonzero proper $\whg$-submodule of
$\mathcal M_\varphi$ (proper because $v\notin U(\whg)\cdot N$).
\emph{Irreducibility when $\kappa_0\neq0$:}
Let $0\neq W\subset\mathcal M_\varphi$ with $w\in W$ of minimal
$\deg_f(w)$.
If $\deg_f(w)=0$: $w\in U(\whh)v\cong M_{\whh,\varphi}$, so
$W\cap M_{\whh,\varphi}\neq0$.  By Theorem~\ref{thm:A}, $v\in W$,
hence $W=\mathcal M_\varphi$.
If $\deg_f(w)\ge1$: let $y$ be a minimal variable dividing
$F_0=\mathrm{LT}_f(w)$.  If $y=f_{n_0}$, Lemma~\ref{lem:reduce-t}
gives $0\neq X\cdot w\in W$ with $\deg_f(X\cdot w)=\deg_f(w)-1$,
contradiction.  If $y=f_{n_0}^1$, Lemma~\ref{lem:reduce-u} gives
the same contradiction.  So $\deg_f(w)=0$ after all.
\end{proof}
\section{The four-point case: general central character}\label{sec:r1-general}
Set $r=1$, so $p(t)=t(t-\alpha)$, $a_1=-\alpha$, $a_2=1$,
$\alpha\in\Cx^\times$.  The center of $\whh$ has basis
$\{1_0,1_1\}$; write $\kappa_0,\kappa_1\in\Cx$ for the respective
scalar actions (no specialization yet).
\subsection{Mixed bracket at $r=1$}
By the Cox--Im bracket formula~\eqref{eq:heis3}, for $r=1$ one has
\begin{equation}\label{eq:psi-r1}
[b_m^1,b_n]=2n\,\psi_{mn}(a),\quad
\psi_{mn}(a)=c_{mn}^{(1)}\omega_1.
\end{equation}
Lemma~\ref{lem:psi-r1-explicit} below gives the complete explicit
formula: $\psi_{mn}(a)=\delta_{m+n,0}\,\omega_1$ for \emph{all} $m,n$
(so $c_{mn}^{(1)}=\delta_{m+n,0}$, $c_{mn}^{(0)}=0$).  In particular,
setting $m=-n$:
\begin{equation}\label{eq:mixed-key}
[b_{-n}^1,b_n]=2n\,\omega_1\quad(n\neq0),
\end{equation}
with \emph{no} $\omega_0$ contribution.
\begin{remark}[Derivation of~\eqref{eq:mixed-key}]\label{rem:mixed-key-proof}
Apply the universal central extension bracket~\eqref{eq:uce-bracket}
with $f=t^{-n}u$ and $g=t^n$.  Since $[h,h]=0$, only the central term
contributes, and $f\,dg = t^{-n}u\cdot nt^{n-1}\,dt = nt^{-1}u\,dt$:
\[
[b_{-n}^1,b_n]
  = (h,h)\,\ol{f\,dg}
  = 2\,\ol{nt^{-1}u\,dt}
  = 2n\,\omega_1,
\]
where we use the basis identification $\omega_1=\ol{t^{-1}u\,dt}$ from
\eqref{eq:omega-def}.  Since $\ol{t^{-1}u\,dt}$ is proportional to
$\omega_1$ only (and not to $\omega_0=\ol{t^{-1}dt}$), the $\omega_0$
component is zero, giving the exact equality above.

Similarly, with $f=t^{n_0}u$ and $g=t^{-n_0}$ (so
$f\,dg = -n_0\,t^{-1}u\,dt$):
\[
[b_{n_0}^1,b_{-n_0}]
  = (h,h)\,\ol{f\,dg}
  = 2\,\ol{(-n_0)t^{-1}u\,dt}
  = -2n_0\,\omega_1,
\]
so $2(-n_0)\psi_{n_0,-n_0}(a)=-2n_0\,\omega_1$ gives
$\psi_{n_0,-n_0}(a)=\omega_1$, confirming $c_{n_0,-n_0}^{(1)}=1$.
\end{remark}
\begin{lemma}[Explicit $\psi_{mn}(a)$ for $r=1$]\label{lem:psi-r1-explicit}
Let $r=1$, $p(t)=t(t-\alpha)$, $a_1=-\alpha$, $a_2=1$.  For all
$m,n\in\Z$:
\begin{equation}\label{eq:psi-explicit-r1}
\psi_{mn}(a) = \delta_{m+n,\,0}\cdot\omega_1,
\end{equation}
equivalently $c_{mn}^{(1)}=\delta_{m+n,0}$ (and $c_{mn}^{(j)}=0$ for
all $j\ge2$, which is vacuous for $r=1$).
In particular, $[b_m^1,b_n]=2n\,\delta_{m+n,0}\,\omega_1$ and
$[b_m,b_n^1]=-2m\,\delta_{m+n,0}\,\omega_1$.
\end{lemma}
\begin{proof}
By the UCE bracket formula~\eqref{eq:uce-bracket} with
$f=t^mu$ and $g=t^n$:
\[
[b_m^1,b_n]
  = (h,h)\,\ol{t^mu\cdot d(t^n)}
  = 2n\,\ol{t^{m+n-1}u\,dt}.
\]
We claim $\ol{t^j u\,dt}=\delta_{j,-1}\,\omega_1$ in $\Omega_R^1/dR$.
For $j\neq-1$: set $k=j+1\neq0$.  Since $t^k u\in R$,
$d(t^k u)=kt^{k-1}u\,dt+t^k\,du$ is exact, hence $\equiv0$
in $\Omega_R^1/dR$.  Since $u\in R$, $du = d(u)$ is also exact,
so $[du]=0$.  Thus $k\,\ol{t^{k-1}u\,dt}=-t^k\,[du]=0$,
giving $\ol{t^j u\,dt}=0$ for all $j\neq-1$.
For $j=-1$: $\ol{t^{-1}u\,dt}=\omega_1\neq0$ by definition of the basis.
Substituting: $[b_m^1,b_n]=2n\,\delta_{m+n-1,-1}\,\omega_1=2n\,\delta_{m+n,0}\,\omega_1$,
so $2n\,\psi_{mn}(a)=2n\,\delta_{m+n,0}\,\omega_1$.
For $n\neq0$ divide by $2n$; for $n=0$ both sides vanish.  This gives
$\psi_{mn}(a)=\delta_{m+n,0}\,\omega_1$.  The formula for $[b_m,b_n^1]$
follows by skew-symmetry.
\end{proof}
\subsection{Extended action formulas for $r=1$}
\begin{lemma}\label{lem:key-action-r1}
For $r=1$, general $(\kappa_0,\kappa_1)$, $m\in S_\varphi^+$,
$n\in S_\varphi^-$, $\ell\ge1$:
\begin{align}
b_m\cdot(b_n)^\ell v
  &=2\ell n\,\kappa_0\,\delta_{m+n,0}\,(b_n)^{\ell-1}v,
  \label{eq:bm-bn-r1}\\
b_m\cdot(b_n^1)^\ell v
  &=-2\ell m\,\delta_{m+n,0}\,\kappa_1\,(b_n^1)^{\ell-1}v,
  \label{eq:bm-bn1-r1}\\
b_m^1\cdot(b_n)^\ell v
  &=2\ell n\,\delta_{m+n,0}\,\kappa_1\,(b_n)^{\ell-1}v,
  \label{eq:bm1-bn-r1}\\
b_m^1\cdot(b_n^1)^\ell v
  &=\ell(n-m)a_{-(m+n)}\kappa_0\,(b_n^1)^{\ell-1}v.
  \label{eq:bm1-bn1-r1}
\end{align}
\end{lemma}
\begin{proof}
As in Lemma~\ref{lem:key-action}, $b_m\cdot v=b_m^1\cdot v=0$ for
$m\in S_\varphi^+$, so $x\cdot y^\ell v=[x,y^\ell]v
=\ell y^{\ell-1}[x,y]v$.  Formula~\eqref{eq:bm-bn-r1} is
identical to~\eqref{eq:bm-bn}.

By Lemma~\ref{lem:psi-r1-explicit}, $\psi_{mn}(a)=\delta_{m+n,0}\,\omega_1$
for all $m,n$ (and $c_{mn}^{(0)}=0$).

For~\eqref{eq:bm-bn1-r1}:
$[b_m,b_n^1]=-[b_n^1,b_m]=-2m\,\psi_{nm}(a)=-2m\,\delta_{m+n,0}\,\omega_1$,
acting on $V$ as $-2m\,\delta_{m+n,0}\,\kappa_1$.

For~\eqref{eq:bm1-bn-r1}:
$[b_m^1,b_n]=2n\,\psi_{mn}(a)=2n\,\delta_{m+n,0}\,\omega_1$,
acting on $V$ as $2n\,\delta_{m+n,0}\,\kappa_1$.

Formula~\eqref{eq:bm1-bn1-r1} is unchanged from~\eqref{eq:bm1-bn1} since
$[b_m^1,b_n^1]=(n-m)a_{-(m+n)}\omega_0$ involves only $\kappa_0$.
\end{proof}
\subsection{Irreducibility locus for $r=1$}
\begin{theorem}\label{thm:A-prime}
Let $r=1$ and assume $\varphi$ is $p$-admissible.  Then
$M_{\whh,\varphi}$ is irreducible if and only if
$(\kappa_0,\kappa_1)\neq(0,0)$.
\end{theorem}
\begin{proof}
\emph{Reducibility when $(\kappa_0,\kappa_1)=(0,0)$:}
All commutators act trivially, so $\whh_\varphi^+$ kills all of
$M_{\whh,\varphi}$ and $M_{\whh,\varphi}^{\ge1}$ is a proper submodule.
\emph{Irreducibility when $(\kappa_0,\kappa_1)\neq(0,0)$:}
Let $0\neq W\subset M_{\whh,\varphi}$ and pick $0\neq w\in W$ with
minimal $\deg_{\whh}(w)=d\ge1$.
As in Theorem~\ref{thm:A}, pick the maximal index $n_0\in S_\varphi^-$
with $\alpha_{n_0}+\beta_{n_0}>0$ in the leading monomial.
\textbf{Case 1: $\alpha_{n_0}>0$.}  Set $m_0:=-n_0\in S_\varphi^+$.
By~\eqref{eq:bm-bn-r1}, $b_{m_0}\cdot w$ has a leading-symbol contribution
$2n_0\kappa_0\cdot\partial_{b_{n_0}}\sigma(w)\neq0$ in
$\gr M_{\whh,\varphi}$ (same associated-graded argument as in
Theorem~\ref{thm:A}).  If $\kappa_0\neq0$, this is a contradiction.
If $\kappa_0=0$ (so $\kappa_1\neq0$): apply $b_{m_0}^1$ instead.
By~\eqref{eq:bm1-bn-r1} with $m_0+n_0=0$ (so $\delta_{m_0+n_0,0}=1$):
$b_{m_0}^1\cdot(b_{n_0})^{\alpha_{n_0}}v
=2\alpha_{n_0}n_0\kappa_1\cdot(b_{n_0})^{\alpha_{n_0}-1}v\neq0$
(since $n_0\in S_\varphi^-$ with $n_0\neq0$, $\kappa_1\neq0$).
The same associated-graded argument shows $b_{m_0}^1\cdot w\neq0$
has degree $d-1$, contradiction.
\textbf{Case 2: $\alpha_{n_0}=0$, $\beta_{n_0}>0$.}  By
$p$-admissibility choose $m_0\in S_\varphi^+$ with $a_{-(m_0+n_0)}\neq0$.
By~\eqref{eq:bm1-bn1-r1}, $b_{m_0}^1\cdot(b_{n_0}^1)^{\beta_{n_0}}v
=\beta_{n_0}(n_0-m_0)a_{-(m_0+n_0)}\kappa_0(b_{n_0}^1)^{\beta_{n_0}-1}v$.
If $\kappa_0\neq0$: nonzero; by associated graded $b_{m_0}^1\cdot w\neq0$
has degree $d-1$, contradiction.
If $\kappa_0=0$ and $\kappa_1\neq0$ in Case~2, set $m_0^*:=-n_0\in S_\varphi^+$.
By~\eqref{eq:bm-bn1-r1} with $m_0^*=-n_0$ (so $\delta_{m_0^*+n_0,0}=1$):
$b_{m_0^*}\cdot(b_{n_0}^1)^{\beta_{n_0}}v
=-2\beta_{n_0}(-n_0)\kappa_1(b_{n_0}^1)^{\beta_{n_0}-1}v
=2\beta_{n_0}n_0\kappa_1(b_{n_0}^1)^{\beta_{n_0}-1}v\neq0$
(since $n_0\neq0$, $\kappa_1\neq0$).
By associated graded $b_{m_0^*}\cdot w\neq0$ has degree $d-1$, contradiction.
In all subcases we reach $\deg_{\whh}(w)=0$, hence $v\in W$.
\end{proof}
\section{Level-zero structure}\label{sec:level-zero}
We now give a complete description of $M_{\whh,\varphi}$ when
$\kappa_0=0$ (and $\kappa_k=0$ for all $k$).
\begin{proposition}\label{prop:level-zero-submodule}
Assume $\kappa_0=\kappa_1=\cdots=\kappa_r=0$.  Then:
\begin{enumerate}[label=(\alph*)]
\item $M_{\whh,\varphi}\cong S(\whh_\varphi^-)\otimes_\Cx V$ as
  $\Cx$-vector spaces, with $\whh_\varphi^-$ acting by polynomial
  multiplication, $\whh_\varphi^+$ acting by zero, and $\whh_0$
  acting by scalars $\lambda,\mu,0,\dots,0$.
\item The filtration $\{F_dM_{\whh,\varphi}\}$ is $\whh$-stable, and
  $M_{\whh,\varphi}^{\ge1}:=F_1M_{\whh,\varphi}
  =\Span\{b^\alpha(b^1)^\beta v:\sum_n\alpha_n+\sum_n\beta_n\ge1\}$
  is a proper nonzero $\whh$-submodule.
\item $M_{\whh,\varphi}^{\ge1}$ corresponds to the augmentation ideal
  $\mathfrak{m}:=\langle b_n,b_n^1:n\in S_\varphi^-\rangle$ of the
  polynomial ring $S(\whh_\varphi^-)$.
\item $M_{\whh,\varphi}^{\ge1}$ is the unique maximal proper submodule
  of $M_{\whh,\varphi}$, and the quotient
  $M_{\whh,\varphi}/M_{\whh,\varphi}^{\ge1}\cong V=\Cx v$ is the unique
  simple quotient.
\item $U(\whg)\cdot M_{\whh,\varphi}^{\ge1}$ is a canonical proper
  nonzero $\whg$-submodule of $\mathcal M_\varphi$.
\end{enumerate}
\end{proposition}
\begin{proof}
(a) Since all $\kappa_k=0$, every commutator in $\whh$ acts as $0$ on any
module (all central elements act trivially).  Hence $\whh_\varphi^+$
annihilates every PBW basis element, and $\whh_\varphi^-$ acts by free
left multiplication: $b_n\cdot(b^\alpha(b^1)^\beta v)=b_n b^\alpha(b^1)^\beta v$
with no correction terms.  The PBW basis
$\{b^\alpha(b^1)^\beta v\}\leftrightarrow\{b^\alpha(b^1)^\beta\otimes v\}$
gives the vector-space isomorphism
$M_{\whh,\varphi}\cong S(\whh_\varphi^-)\otimes_\Cx V$.
(b) Since $\whh_\varphi^+$ acts as $0$, and $\whh_\varphi^-$ maps
$F_dM$ to $F_{d+1}M$ (each generator application increases PBW degree
by $1$), the filtration $\{F_dM_{\whh,\varphi}\}$ is $\whh$-stable.
$M_{\whh,\varphi}^{\ge1}=F_1M_{\whh,\varphi}$ is nonzero (it contains
$b_n v$ for any $n\in S_\varphi^-$) and proper (it does not contain $v$,
which has PBW degree $0$).
(c) Under the identification of (a),
$M_{\whh,\varphi}^{\ge1}=F_1M_{\whh,\varphi}$ corresponds to
$\mathfrak{m}\otimes_\Cx V$ where $\mathfrak{m}=\langle b_n,b_n^1:
n\in S_\varphi^-\rangle$ is the augmentation ideal of $S(\whh_\varphi^-)$.
(d) The PBW-degree-$0$ piece $F_0M_{\whh,\varphi}=\Cx v$ is
one-dimensional.  Any submodule $W\not\ni v$ satisfies
$W\cap F_0M=0$, hence $W\subset F_1M=M_{\whh,\varphi}^{\ge1}$.
Since $M_{\whh,\varphi}^{\ge1}$ is itself such a submodule, it is
the unique maximal one.
(e) Write $N:=U(\whg)\cdot M_{\whh,\varphi}^{\ge1}$.
We claim $v\notin N$.  By the $f$-PBW decomposition
(Lemma~\ref{lem:pbw-Mphi}), every element of $N$ can be written as
$\sum_F Fw_F$ with $w_F\in U(\whh)\cdot M_{\whh,\varphi}^{\ge1}$.
If $\deg_f(w)=0$ for some nonzero $w\in N$, then
$w\in U(\whh)\cdot M_{\whh,\varphi}^{\ge1}\subset F_1M_{\whh,\varphi}$
(since $M_{\whh,\varphi}^{\ge1}$ is $\whh$-stable and every element
has PBW degree $\ge1$).  In particular $w\neq v$, so $v\notin N$.
Hence $N\subsetneq\mathcal M_\varphi$.
\end{proof}
\section{Examples and admissibility}\label{sec:examples}
\subsection{An explicit $p$-admissible polarization}
\begin{example}[All-negative polarization]\label{ex:all-minus}
Let $\varphi\equiv-$, so $\varphi(n)=-$ for all $n>0$.
Then $S_\varphi^-=\{1,2,3,\dots\}$ and $S_\varphi^+=\{-1,-2,-3,\dots\}$.
By Remark~\ref{rem:standard-admissible}, $\varphi\equiv-$ is
$p$-admissible: for any $n\in S_\varphi^-$, choose $m:=-n-(r+1)\in S_\varphi^+$
and verify $a_{-(m+n)}=a_{r+1}=1\neq0$.
Theorems~\ref{thm:A} and~\ref{thm:B} therefore apply.
\end{example}
\subsection{The four-point case $r=1$}
\begin{example}[$r=1$]\label{ex:r=1}
$p(t)=t^2-\alpha t$, $a_1=-\alpha$, $a_2=1$.
The commutator $[b_m^1,b_n^1]$ is supported on $m+n\in\{-1,-2\}$.
Admissibility requires: for each $n\in S_\varphi^-$, some $m\in S_\varphi^+$
with $m+n\in\{-1,-2\}$.
\end{example}
We now give a worked example illustrating irreducibility directly.
\begin{example}[Explicit computation, $r=1$, $\varphi\equiv-$]\label{ex:worked}
Let $r=1$, $\alpha=1$ (so $p(t)=t^2-t$, $a_1=-1$, $a_2=1$),
polarization $\varphi\equiv-$ (so $S_\varphi^-=\{1,2,\dots\}$,
$S_\varphi^+=\{-1,-2,\dots\}$), $\lambda=\mu=0$, $\kappa_0=1$.
The degree-$1$ piece of $M_{\whh,\varphi}$ has basis $\{b_1v,\,b_1^1v\}$.
We compute the action of the degree-$(-1)$ operators
$b_{-1},b_{-1}^1\in\whh_\varphi^+$ on these basis vectors.
\textit{Action on $b_1 v$:}
\begin{itemize}
\item From~\eqref{eq:bm-bn}: $b_{-1}\cdot(b_1 v)=2(1)\kappa_0\delta_{0,0}\,v
  =2v\neq0$.
\item From~\eqref{eq:bm1-bn}: $b_{-1}^1\cdot(b_1 v)=0$ (mixed bracket
  lies in $\Span\{1_1\}$, acting as $0$ since $\kappa_1=0$).
\end{itemize}
\textit{Action on $b_1^1 v$:}
\begin{itemize}
\item From~\eqref{eq:bm-bn1}: $b_{-1}\cdot(b_1^1 v)=0$.
\item From~\eqref{eq:bm1-bn1}: $b_{-1}^1\cdot(b_1^1 v)
  =(1-(-1))\,a_{-((-1)+1)}\,\kappa_0\,v
  =2\cdot a_0\cdot 1\cdot v=0$ (since $a_0=0$).
  We need $m+n=-1$: choose $m=-2\in S_\varphi^+$:
  $b_{-2}^1\cdot(b_1^1 v)=(1-(-2))\,a_{-((-2)+1)}\,\kappa_0\,v
  =3\cdot a_{1}\cdot 1\cdot v=3(-1)v=-3v\neq0$.
\end{itemize}
Thus $b_{-1}\cdot(b_1 v)=2v\neq0$: any submodule containing $b_1 v$
contains $v$, hence equals $M_{\whh,\varphi}$.
Likewise $b_{-2}^1\cdot(b_1^1 v)=-3v\neq0$: any submodule containing
$b_1^1 v$ contains $v$.  Every degree-$1$ basis vector generates the
whole module, confirming irreducibility at $\kappa_0=1\neq0$.
\end{example}
\subsection{Non-admissible polarizations: the sharp example}
\begin{example}[Non-admissible]\label{ex:nonadmissible}
Let $r=1$ and define $\varphi$ by $\varphi(1)=-$ and $\varphi(n)=+$
for all $n>1$ (extended by oddness).  Then
$S_\varphi^-=\{1,-2,-3,\dots\}$ and $S_\varphi^+=\{-1,2,3,\dots\}$.
The index $n=1\in S_\varphi^-$ requires some $m\in S_\varphi^+$ with
$a_{-(m+1)}\neq0$, i.e.\ $m+1\in\{1,2\}$, i.e.\ $m\in\{0,1\}$.
But $0\notin S_\varphi^+$ and $1\notin S_\varphi^+$ (since $\varphi(1)=-$
gives $1\in S_\varphi^-$, not $S_\varphi^+$).  Thus $\varphi$ is
\emph{not} $p$-admissible.
\end{example}
\begin{proposition}\label{prop:nonadmissible-reducible}
For the non-admissible $\varphi$ of Example~\ref{ex:nonadmissible},
with $r=1$, there exist $\kappa_0\neq0$ and nonzero $(\lambda,\mu)$ such
that $M_{\whh,\varphi}$ is reducible.
\end{proposition}
\begin{proof}
We work under $\kappa_1=0$ and show that $W:=\Cx b_1^1 v$ is
$\whh$-stable, hence a proper submodule not containing $v$.
Recall $S_\varphi^+=\{-1,2,3,\dots\}$ for this $\varphi$.
We must verify $x\cdot(b_1^1v)=0$ for all generators $x\in\whh_\varphi^+$.
\emph{Type $b_m^1$ ($m\in S_\varphi^+$)}:
$b_m^1\cdot(b_1^1v)=(1-m)a_{-(m+1)}\kappa_0\,v$ by~\eqref{eq:bm1-bn1}.
We need $a_{-(m+1)}\neq0$, i.e.\ $-(m+1)\in\{1,2\}$, i.e.\
$m\in\{-2,-3\}$.  Checking: $-2\notin S_\varphi^+=\{-1,2,3,\dots\}$ and
$-3\notin S_\varphi^+$.  For all other $m\in S_\varphi^+$:
$a_{-(m+1)}=0$, so $b_m^1\cdot(b_1^1 v)=0$. \checkmark
\emph{Type $b_m$ ($m\in S_\varphi^+$)}: by~\eqref{eq:bm-bn1},
$b_m\cdot(b_1^1 v)$ picks up $-2m\,\psi_{1,m}(a)$ acting on $v$.
By Lemma~\ref{lem:psi-r1-explicit}, $\psi_{1,m}(a)=\delta_{1+m,0}\,\omega_1$
for all $m$.  For $m\neq-1$: $\delta_{1+m,0}=0$, so the bracket acts
as $0$.  For $m=-1\in S_\varphi^+$: $\psi_{1,-1}(a)=\omega_1$, and
$[b_{-1},b_1^1]=-2(-1)\psi_{1,-1}(a)=2\omega_1$, acting as $2\kappa_1=0$.
Hence $b_m\cdot(b_1^1 v)=0$ for all $m\in S_\varphi^+$. \checkmark
Therefore $\whh_\varphi^+\cdot(b_1^1 v)=0$ for any $\kappa_0\neq0$.
Let $W:=U(\whh)\cdot b_1^1v$ be the $\whh$-submodule generated by $b_1^1v$.
Since $\whh_\varphi^+$ kills $b_1^1v$ and every commutator in $\whh$ is
central, the element $b_1^1v$ has no degree-lowering contributions: any
operator in $\whh_\varphi^+$ applied to $U(\whh)\cdot b_1^1v$ either kills
$b_1^1v$ directly (as shown) or, when commuted through to the right, produces
central scalars times lower-degree terms in $W$, all of PBW degree $\ge1$.
Precisely: the PBW filtration is $\whh$-stable, $b_1^1v\in F_1M$, and
$\whh_\varphi^+\cdot F_1M\subset F_0M=\Cx v$ requires producing $v$; but
we showed $\whh_\varphi^+\cdot b_1^1v=0$, so in particular $W\subset F_1M$
and $v\notin W$.  Hence $W$ is a proper nonzero $\whh$-submodule, so
$M_{\whh,\varphi}$ is reducible, confirming that $p$-admissibility is
necessary for Theorem~\ref{thm:A}.
\end{proof}
\section{Higher-dimensional tops}\label{sec:higher-tops}
We show that the irreducibility criterion extends to finite-dimensional
tops.
\begin{definition}\label{def:top-fd}
A \emph{finite-dimensional top} is a finite-dimensional
$\wh{\mathfrak b}_\varphi$-module $\mathcal{V}$ on which
$\whh_\varphi^+$ acts trivially and $\whh_0$ acts by scalars
(i.e.\ $\mathcal{V}$ is a direct sum of one-dimensional top modules
$V_j=\Cx v_j$ of Definition~\ref{def:top1d}).
Set $M_{\whh,\varphi}(\mathcal{V}):=U(\whh)\otimes_{U(\wh{\mathfrak
b}_\varphi)}\mathcal{V}$.
\end{definition}
\begin{corollary}\label{cor:higher-top}
Assume $\kappa_1=\cdots=\kappa_r=0$ and $\varphi$ is $p$-admissible.
Let $\mathcal{V}=\bigoplus_j V_j$ be a finite-dimensional top.  Then:
\begin{enumerate}[label=(\alph*)]
\item If $\kappa_0\neq0$, then $M_{\whh,\varphi}(\mathcal{V})$ is
  irreducible if and only if $\mathcal{V}$ is an irreducible
  $\whh_0$-module.
\item If $\kappa_0=0$, then $M_{\whh,\varphi}(\mathcal{V})$ is reducible.
\end{enumerate}
\end{corollary}
\begin{proof}
(a) \emph{Reducibility of $\mathcal V$ $\Rightarrow$ reducibility of
$M_{\whh,\varphi}(\mathcal V)$:}
Suppose $\mathcal V=V_1\oplus V_2$ with both summands nonzero.  By the PBW
theorem applied to $M_{\whh,\varphi}(\mathcal V)=U(\whh)\otimes_{U(\wh{\mathfrak b}_\varphi)}\mathcal V$,
there is a vector-space isomorphism
$M_{\whh,\varphi}(\mathcal V)\cong S(\whh_\varphi^-)\otimes_\Cx\mathcal V
\cong(S(\whh_\varphi^-)\otimes V_1)\oplus(S(\whh_\varphi^-)\otimes V_2)$.
The summand $S(\whh_\varphi^-)\otimes V_2\cong M_{\whh,\varphi}(V_2)$ is
preserved by $\whh_\varphi^-$ (which maps it to higher PBW degree within
the same summand), by $\whh_0$ (which acts by scalars on $V_2$), and by
$\whh_\varphi^+$ (which kills $V_2$ by definition of top module).  Hence
$M_{\whh,\varphi}(V_2)$ is a nonzero proper $\whh$-submodule.

\emph{Irreducibility of $\mathcal V$ $\Rightarrow$ irreducibility of
$M_{\whh,\varphi}(\mathcal V)$:}
Since $\whh_0$ is abelian, every irreducible finite-dimensional
$\whh_0$-module is one-dimensional, so $\mathcal V=V=\Cx v$.
Theorem~\ref{thm:A} applies directly: $\kappa_0\neq0$ and $p$-admissibility
imply $M_{\whh,\varphi}(V)$ is irreducible.

(b) Since all $\kappa_k=0$, every commutator acts as $0$; by the same
argument as Theorem~\ref{thm:A}, $F_1M_{\whh,\varphi}(\mathcal{V})$ is a
proper $\whh$-submodule.
\end{proof}
\section{Outlook and further directions}\label{sec:outlook}
\begin{itemize}[leftmargin=2em]
\item \textbf{General central characters for $r\ge2$.}
The four-point analysis of Section~\ref{sec:r1-general} extends to
general $r$: under $p$-admissibility, $M_{\whh,\varphi}$ should be
irreducible iff $(\kappa_0,\kappa_1,\dots,\kappa_r)\neq\mathbf{0}$.
Making this precise requires tracking the full Cox--Im formulas for
$\psi_{mn}(a)$ and establishing appropriate generalizations of
Lemma~\ref{lem:key-action-r1}.
\item \textbf{Other simple Lie algebras.}
The $\sltwo$ framework (and the explicit Cox--Im presentation) could be
extended to $\g$ simple: the effective level and the admissibility
condition would involve the Killing form of $\g$ and the root system
structure.  The KN geometry should govern simplicity in analogous ways.
\item \textbf{Elliptic and DJKM settings.}
The $p$-admissibility condition and the nonzero-level mechanism should
have natural analogues in the elliptic~\cite{Bueno:2009aa} and
DJKM~\cite{date1983landau,Cox:2013ab,cox2013djkm} settings, where similar Heisenberg
subalgebras and imaginary mode structures appear.
\item \textbf{Category $\mathcal{O}$ and Verma module comparisons.}
It would be natural to embed the $\varphi$-Verma modules into a suitable
BGG-type category $\mathcal{O}$ for $\whg$ and compare the
$\varphi$-parameterized family with the classical Verma modules and
their composition series.
\end{itemize}
\subsection*{Acknowledgements}
Part of this work was discussed with Vyacheslav Futorny. I am grateful to
him for valuable comments and suggestions.  This study was financed, in
part, by the S\~ao Paulo Research Foundation (FAPESP), grant 2024/14914-9.

\end{document}